\input amstex
\documentstyle{amsppt}
\output={\plainoutput}
\pageno=1
\footline={\ifnum\pageno>1 \myfl\fi}
\define\myfl{\hss\tenrm\folio\hss}
\headline={\hfil}
\NoBlackBoxes
\magnification=\magstep 1
\pagewidth{15 true cm}
\pageheight{22 true cm}

\topmatter
\title The ratio of two zeta-determinants of Dirac Laplacians associated with unitary involutions on a compact manifold
with cylindrical end
\endtitle
\author Yoonweon Lee
\endauthor
\affil  Department of Mathematics \\
        Inha University \\
        Incheon, 402-751, Korea
\endaffil
\endsubjclass
\thanks
This work was supported by KRF-2003-041-C00028.
\endthanks
\abstract
Given two unitary involutions $\sigma_{1}$ and $\sigma_{2}$ satisfying $G \sigma_{i} = - \sigma_{i} G$
on $ker B$
on a compact manifold with cylindrical end, M. Lesch, K. Wojciechowski ([LW]) and W. M\"uller ([M]) established the formula
describing the difference of two eta-invariants
with the APS boundary conditions associated with $\sigma_{1}$ and $\sigma_{2}$. In this paper we establish the analogous formula
for the zeta-determinants of Dirac Laplacians. For the proof of the result we use the Burghelea-Friedlander-Kappeler's
gluing formula for zeta-determinants and the scattering theory developed by W. M\"uller in [M].
This result was also obtained independently by J. Park and K. Wojciechowski ([PW2]).
\endabstract
\endtopmatter

\document
\baselineskip 0.6 true cm plus 3 pt minus 3 pt

\S 1 {\bf Introduction}
\TagsOnRight
\vskip 0.5 true cm

The purpose of this paper is to discuss the ratio of two zeta-determinants of Dirac Laplacians with the
Atiyah-Patodi-Singer (APS) boundary conditions associated with two unitary involutions.
This subject has already been studied by J. Park and K. P. Wojciechowski in [PW1] and [PW2]. However, we here present a
completely different way of proving the same result by using the Burghelea-Friedlander-Kappeler's gluing formula
(BFK-gluing formula) for zeta-determinants developed in [L2] and [L3]. The motivation for this analysis comes from the works of
M. Lesch, K. Wojciechowski ([LW]) and W. M\"uller ([M]).
In these papers they established independently the formula describing the difference of two eta-invariants by mod ${\Bbb Z}$
with the APS boundary conditions associated with unitary involutions $\sigma_{i}$ ($i = 1, 2$) on $ker B$ satisfying $\sigma_{i}G = - G\sigma_{i}$
on a compact manifold with cylindrical end.

In this paper we discuss the analogous question for the zeta-determinants of Dirac Laplacians. More precisely we study the ratio of two
zeta-determinants of Dirac Laplacians with the APS boundary conditions associated with $\sigma_{1}$ and $\sigma_{2}$.
As the arguments in [PW2] (or [PW1]), the proof of the main result consists of two parts.
We first investigate the zeta-determinants of the so-called Dirichlet-to-Neumann operators appearing in the BFK-gluing formula
by using the scattering theory in [M] and then prove that the ratio
of two zeta-determinants does not depend on the cylinder length.
Next, we use the adiabatic limits with respect to the cylinder length to obtain the main result, again with the help of the scattering theory in [M].

\vskip 0.3 true cm

Now we introduce the basic settings.
Let $(M, g)$ be a compact oriented $m$-dimensional Riemannian manifold ($m>1$) with boundary $Y$
and $E\rightarrow M$ be a Clifford module bundle.
Choose a collar neighborhood $N$ of $Y$ which is diffeomorphic to $(-1,0] \times Y$.
We assume that the metric $g$ is a product one on $N$ and the bundle $E$ has the product structure
on $N$, {\it i.e.} $E|_{N}=p^{\ast}E|_{Y}$, where $p : (-1, 0] \times Y \rightarrow Y$ is the canonical projection.
Suppose that $D_{M}$ is a compatible Dirac operator acting on smooth sections of $E$.
We assume that $D_{M}$ has the following form on $N$.
$$  D_{M}= G( \partial_{u}+B),$$
where $G : E|_{Y} \rightarrow E|_{Y}$ is a bundle automorphism, $\partial_{u}$ is the outward normal derivative to $Y$
on $N$ and $B$ is a Dirac operator on $Y$.
We further assume that $G$ and $B$ are independent of the normal coordinate $u$ and satisfy
$$G^{\ast}=-G, \qquad G^{2}=-I, \qquad B^{\ast}=B  \qquad {\text and }  \qquad GB=-BG.$$
Then we have, on $N$,
$$ D_{M}^{2}=-\partial_{u}^{2}+B^{2}.$$

We now define the APS boundary condition $P_{<}$ (or $P_{>}$) by the orthogonal projection
to the space spanned by negative (or positive) eigensections of $B$.
If $ker B \neq \{ 0 \}$, we need an extra condition on $P_{<}$, say, a unitary involution on $ker B$ anticommuting with $G$.
Suppose that $\sigma : ker B \rightarrow ker B$ is a unitary involution satisfying
$$\sigma G = - G \sigma,  \qquad \qquad  \sigma^{2} = I_{Ker B}.$$
We denote by $P_{\sigma}$
$$ P_{\sigma}= P_{<}  + \frac{1}{2} (I - \sigma)|_{ker B}.$$
Then $D_{M, P_{\sigma}}$, the Dirac operator $D_{M}$ with the domain
$Dom{D_{M, P_{\sigma}}} = \{ \phi \in C^{\infty}(E) \mid  P_{\sigma} (\phi|_{Y})=0 \} $,
is an essential self-adjoint elliptic operator having discrete real spectrum.
Denoting by $\eta_{D_{M, P_{\sigma}}}(0)$ the eta-invariant associated with $D_{M,P_{\sigma}}$, we have
the following theorem, which is due to Lesch, Wojciechowski ([LW]) and M\"uller ([M]).

\proclaim{Theorem 1.1}
Suppose that $\sigma_{1}$, $\sigma_{2} : ker B \rightarrow ker B$ are two unitary involutions
satisfying $G \sigma_{i} = - \sigma_{i} G$ on $ker B$ and
put  $ P_{\sigma_{i}}= P_{<} + \frac{1}{2} (I - \sigma_{i})|_{ker B}$  ($i= 1, 2$).
Then
$$\eta_{D_{M, P_{\sigma_{1}}}}(0) - \eta_{D_{M, P_{\sigma_{2}}}}(0) \equiv
-\frac{1}{\pi i} \log det (\sigma_{2}\sigma_{1}|_{ker(G-i)})
\qquad  (mod {\Bbb Z}) .$$
\endproclaim

\vskip 0.3 true cm

In this paper we are going to establish the analogous result of Theorem 1.1 for the zeta-determinants of Dirac Laplacians
under the assumption that both $D_{M, P_{\sigma_{1}}}$ and $D_{M, P_{\sigma_{2}}}$  are invertible.
The Dirac Laplacian $D_{M,P_{\sigma_{i}}}^{2}$ is defined to be the operator $D_{M}^{2}$ with the following domain.
$$Dom(D_{M, P_{\sigma_{i}}}^{2}) =
\{ \phi \in C^{\infty}(E) \mid  P_{\sigma_{i}} (\phi|_{Y})=0, P_{\sigma_{i}}((D_{M}\phi)|_{Y}) = 0  \}.$$
Then we are going to discuss how to describe
$$ \frac{Det D_{M,P_{\sigma_{1}}}^{2}}{Det D_{M, P_{\sigma_{2}}}^{2}}.  \tag1.1 $$
For this purpose we are going to use the BFK-gluing formula as a basic tool.
Since the BFK-gluing formula works best for invertible operators, we assume that
both $D_{M, P_{\sigma_{1}}}$ and $D_{M, P_{\sigma_{2}}}$  are invertible operators.

Let us introduce the basic settings of the adiabatic limit.
We denote by $M_{r}$ the manifold with boundary obtained by attaching $[0, r] \times Y$ to $M$,
where we identify $Y$ with $Y_{0} := \{ 0 \} \times Y$,
{\it i.e.} $M_{r} = M\cup_{Y} [0, r] \times Y$. The bundle $E\rightarrow M$ and the Dirac operator $D_{M}$
can be extended naturally to the bundle $E_{r}\rightarrow M_{r}$ and the Dirac operator $D_{M_{r}}$ on $M_{r}$ by the product
structures. We also denote by $M_{\infty} := M\cup_{Y} [0, \infty)\times Y$ and by $D_{M_{\infty}}$ the extension of $D_{M}$ to $M_{\infty}$.
Then the scattering theory in [M] shows that
the Dirac operator $D_{M_{\infty}}$ determines the
unitary involution $C(0)$ anticommuting with $G$ on $ker B$, which is called the scattering matrix.
The following theorems are the main results of this paper,
which are proved in Section 4 by using the BFK-gluing formula and the scattering theory.

\proclaim{Theorem 1.2}
Suppose that both $D_{M, P_{\sigma_{1}}}^{2}$ and $D_{M, P_{\sigma_{2}}}^{2}$ are invertible operators. Then
$$\frac{Det D_{M_{r}, P_{\sigma_{1}}}^{2}}{Det D_{M_{r}, P_{\sigma_{2}}}^{2}}$$
does not depend on the cylinder length $r$.
\endproclaim

\noindent
By taking the limit of $\frac{DetD_{M_{r}, P_{\sigma_{1}}}^{2}}{DetD_{M_{r}, P_{\sigma_{2}}}^{2}}$ as $r\rightarrow \infty$,
we have the following theorem.

\proclaim{Theorem 1.3}
Suppose that both $D_{M, P_{\sigma_{1}}}^{2}$ and $D_{M, P_{\sigma_{2}}}^{2}$ are invertible operators.
Then
$$
\frac{ Det D_{M, P_{\sigma_{1}}}^{2} }{ Det D_{M, P_{\sigma_{2}}}^{2} } =
 \frac{det \left( C(0) - \sigma_{1} \right) }{det \left( C(0) - \sigma_{2} \right)} ,
$$
where $det \left( C(0) - \sigma_{i} \right)$ is the usual determinant of a linear map $ C(0) - \sigma_{i} $ on $ker B$
and $C(0)$ is the scattering matrix explained above.
\endproclaim

\vskip 0.3 true cm
\noindent
{\it Remark} : \hskip 0.1 true cm
(1) The invertibility of $D_{M, P_{\sigma_{i}}}^{2}$  ($D_{M, P_{\sigma_{i}}}$) is equivalent to
the non-existence of the $L^{2}$-solutions of $D_{M_{\infty}}$
on $M_{\infty}$ and the invertibility of $ C(0) - \sigma_{i} $. In this case one can check easily that
$ker ( I - \sigma_{i} ) \cap ker ( I - C(0) ) = \{ 0 \} $.  \newline
(2) This result was proved independently by Park and Wojciechowski in [PW2] ({\it cf.} [PW1]).

\vskip 1 true cm

\S 2 {\bf The BFK-gluing formula for zeta-determinants }

\vskip 0.5 true cm

In this section we are going to describe the BFK-gluing formula developed in [L2] and [L3]. Recall that
$M_{r} = M \cup_{Y} [0, r]\times Y$ and define $Q_{1}, Q_{2, r} : C^{\infty}(Y_{0})\rightarrow C^{\infty}(Y_{0})$
as follows.
For any $f\in C^{\infty}(Y_{0})$, choose $\phi_{1}\in C^{\infty}(M)$ and $\phi_{2,r}\in C^{\infty}([0, r]\times Y)$
satisfying
$$D_{M}^{2}\phi_{1}=0,  \quad  (-\partial_{u}^{2}+B^{2})\phi_{2,r}=0,$$
$$\phi_{1}|_{Y_{0}}=\phi_{2,r}|_{Y_{0}}=f,    \quad
P_{\sigma}(\phi_{2,r}|_{Y_{r}}) = P_{\sigma}(G(\partial_{u}+B)(\phi_{2,r})|_{Y_{r}}) = 0,$$
where $P_{\sigma}= P_{>} + \frac{1}{2}(I - \sigma)$ for a given unitary involution $\sigma$ satisfying $ G \sigma = - \sigma G$ on $ker B$
and $Y_{0} := \{ 0 \} \times Y$, $Y_{r} := \{ r \} \times Y$.
Then we define
$$Q_{1}(f) = (\partial_{u}\phi_{1})|_{Y_{0}} , \qquad Q_{2,r}(f)=-(\partial_{u}\phi_{2,r})|_{Y_{0}}$$
and define the Dirichlet-to-Neumann operator $R_{r, \sigma} : C^{\infty}(Y_{0})\rightarrow C^{\infty}(Y_{0})$ by
$$R_{r, \sigma}(f)= Q_{1}(f)-Q_{2,r}(f) = (\partial_{u}\phi_{1})|_{Y_{0}} - (\partial_{u}\phi_{2,r})|_{Y_{0}}.$$
It is known in [L1] that the Dirichlet-to-Neumann operator $R_{r, \sigma}$ is an elliptic $\Psi$DO of order 1 and
if $D^{2}_{M_{r},P_{\sigma}}$
is invertible, $R_{r, \sigma}$ is also invertible.
We denote by $\gamma_{0}$ ($\gamma_{r}$)
the Dirichlet boundary condition on $Y_{0}$ ($Y_{r}$) and also denote by $D^{2}_{M, \gamma_{0}}$ the Dirac Laplacian with the domain
$$Dom(D^{2}_{M, \gamma_{0}}) = \{ \phi \in C^{\infty}(M) \mid \phi|_{Y} = 0 \}.$$
Similarly, we denote by $(-\partial_{u}^{2}+B^{2})_{\gamma_{0}, P_{\sigma}}$ the Dirac Laplacian on $[0, r]\times Y$ with the domain
$$ Dom((-\partial_{u}^{2}+B^{2})_{\gamma_{0}, P_{\sigma}}) = $$
$$ \{ \phi \in C^{\infty}([0, r] \times Y) \mid \phi|_{Y_{0}} = 0,
P_{\sigma}(\phi|_{Y_{r}}) = 0 \text{ and } P_{\sigma}(G(\partial_{u} + B) \phi|_{Y_{r}}) = 0 \}.$$
Then the following theorem is due to Burghelea, Friedlander, Kappeler and the author for the constant part,
which we call the BFK-gluing formula for zeta-determinants
([BFK], [L1], [L2]).

\proclaim{Theorem 2.1}
Suppose that $D^{2}_{M_{r},P_{\sigma}}$ is an invertible Dirac Laplacian on $M_{r}$. Then
$$ \log Det D^{2}_{M_{r},P_{\sigma}} = \log Det D^{2}_{M, \gamma_{0}}
+ \log Det (-\partial_{u}^{2}+B^{2})_{\gamma_{0}, P_{\sigma}} $$
$$ - \log2 \cdot (\zeta_{B^{2}}(0) + dim ker B) + \log Det R_{r, \sigma}.$$
\endproclaim

\vskip 0.3 true cm

The following lemma describes the spectrum of $(-\partial_{u}^{2}+B^{2})_{\gamma_{0}, P_{\sigma}}$, which can be computed straightforwardly.

\proclaim{Lemma 2.2}
The spectrum of $(-\partial_{u}^{2}+B^{2})_{\gamma_{0}, P_{\sigma}}$ is as follows.
$$Spec \left( (-\partial_{u}^{2}+B^{2})_{\gamma_{0}, P_{\sigma}} \right) = $$
$$ \left\{ \lambda^{2} + (\frac{k \pi}{r})^{2} \mid \lambda\in Spec(B), \lambda < 0, k = 1, 2, 3, \cdots \right\}
\cup \left\{ \left( \frac{k \pi}{r} \right)^{2} \mid k = 1, 2, 3, \cdots \right\} $$
$$ \cup \quad \left\{ \left( \frac{(k+\frac{1}{2})\pi}{r} \right)^{2} \mid k = 0, 1, 2, \cdots \right\}  \quad
 \cup \quad \{ \mu_{\lambda, j} \mid \lambda\in Spec(B), \lambda > 0, \mu_{\lambda, j} > \lambda^{2} \},
\qquad \qquad \qquad $$
where $\mu_{\lambda, j}$'s are the solutions of the following equation with $\lambda > 0 $
$$\sqrt{\mu - \lambda^{2}} cos (\sqrt{\mu - \lambda^{2}} r) + \lambda sin(\sqrt{\mu - \lambda^{2}} r) = 0$$
and the multiplicities of $\left(\frac{(k+\frac{1}{2})\pi}{r} \right)^{2}$
and $\left(\frac{k \pi}{r} \right)^{2}$ are $\frac{1}{2}dim ker B$.
\endproclaim

\vskip 0.3 true cm

\noindent
Lemma 2.2 shows that on a cylinder the eigenvalues of $(-\partial_{u}^{2}+B^{2})_{\gamma_{0}, P_{\sigma}}$
are independent of the choice of unitary involutions on $ker B$.

Now suppose that we are given two unitary involutions $\sigma_{1}$ and $\sigma_{2}$ anticommuting with $G$ on $ker B$.
Assume that $D^{2}_{M, P_{\sigma_{1}}}$, $D^{2}_{M, P_{\sigma_{2}}}$ are invertible operators.
Then the {\it Remark} below Theorem 1.3 shows that
$D^{2}_{M_{r}, P_{\sigma_{1}}}$, $D^{2}_{M_{r}, P_{\sigma_{2}}}$ are invertible operators for each $r > 0$.
Hence, Theorem 2.1 and Lemma 2.2 lead to
$$\log Det D^{2}_{M_{r}, P_{\sigma_{1}}} - \log Det D^{2}_{M_{r}, P_{\sigma_{2}}} = \log Det R_{r,\sigma_{1}}
- \log Det R_{r,\sigma_{2}} . \tag2.1 $$
\noindent
To analyze the equation (2.1)  we first investigate the operator $R_{r,\sigma_{i}}$ ($i= 1, 2$).
The following lemma is also straightforward.

\proclaim{Lemma 2.3}
For any $f \in C^{\infty}(Y)$ with $Bf=\lambda f$, $R_{r, \sigma_{i}} (f)$ ($ i = 1, 2$) is given as follows.
$$
R_{r, \sigma_{i}} (f) = \cases
Q_{1}(f) + \left( |\lambda| + \frac{2|\lambda| e^{-r|\lambda|}}{e^{r|\lambda|} - e^{-r|\lambda|}} \right) f & \text{ if $\lambda < 0$} \\
Q_{1}(f) + |\lambda| f  & \text{ if $\lambda > 0$} \\
Q_{1}(f) + \frac{1}{2r}(I - \sigma_{i}) f   & \text{ if $\lambda = 0$}
\endcases
$$
\endproclaim

\noindent
Taking the derivative on $R_{r, \sigma_{i}}$ with respect to $r$ gives the following corollary.

\proclaim{Corollary 2.4} For $i= 1, 2$
$$\frac{d}{dr} R_{r, \sigma_{i}} = -\frac{4|B|^{2}}{\left( e^{r|B|} - e^{-r|B|} \right)^{2}} P_{<} -
\frac{1}{2r^{2}}(I - \sigma_{i})|_{ker B},$$
and hence
$ R_{r, \sigma_{i}}^{-1} (\frac{d}{dr} R_{r, \sigma_{i}}) $ is a trace class operator for each $r > 0$.
\endproclaim

Setting $K_{r} = \frac{4|B|^{2}}{\left( e^{r|B|} - e^{-r|B|} \right)^{2}}$, we have
$$ \frac{d}{dr} \left( \log Det R_{r, \sigma_{1}} - \log Det R_{r, \sigma_{2}} \right)
 = Tr \left( R_{r, \sigma_{1}}^{-1} \frac{d}{dr} R_{r, \sigma_{1}} - R_{r, \sigma_{2}}^{-1} \frac{d}{dr} R_{r, \sigma_{2}} \right) $$
$$
\multline
= \frac{1}{2r^{2}} Tr \left\{ R_{r, \sigma_{2}}^{-1} (I - \sigma_{2})|_{ker B} -
R_{r, \sigma_{1}}^{-1} (I - \sigma_{1})|_{ker B} \right\}  \\
+ Tr \left\{ \left( R_{r, \sigma_{2}}^{-1} - R_{r, \sigma_{1}}^{-1} \right) K_{r} P_{<} \right\} .
\endmultline  \tag2.2
$$
\noindent
We can see easily from Lemma 2.3 that
$$
R_{r, \sigma_{1}} - R_{r, \sigma_{2}} = \frac{1}{2r} (\sigma_{2}-\sigma_{1})|_{ker B}, \qquad \qquad  \tag2.3 $$
and hence we have
$$
\split
Tr \left\{ \left( R_{r, \sigma_{2}}^{-1} - R_{r, \sigma_{1}}^{-1} \right) K_{r} P_{<} \right\}
 & = Tr \left\{  R_{r, \sigma_{2}}^{-1} (R_{r, \sigma_{1}} -
R_{r, \sigma_{2}}) R_{r, \sigma_{1}}^{-1}  K_{r} P_{<} \right\} \\
& = \frac{1}{2r} Tr \left\{ (\sigma_{2} - \sigma_{1})|_{ker B} R_{r, \sigma_{1}}^{-1}
K_{r} P_{<} R_{r, \sigma_{2}}^{-1} \right\}.
\endsplit \tag2.4 $$
\noindent
From (2.1) to (2.4) we have the following lemma.

\proclaim{Lemma 2.5}
$$
\frac{d}{dr} \left\{ \log Det D^{2}_{M_{r}, P_{\sigma_{1}}} - \log Det D^{2}_{M_{r}, P_{\sigma_{2}}} \right\} =  $$
$$ \multline
\frac{1}{r^{2}} Tr \left\{ R_{r, \sigma_{2}}^{-1} proj_{ker(I + \sigma_{2})} -
R_{r, \sigma_{1}}^{-1} proj_{ker(I + \sigma_{1})} \right\} +   \\
\frac{1}{2r} Tr \left\{ (\sigma_{2} - \sigma_{1})|_{ker B} R_{r, \sigma_{1}}^{-1}
K_{r} P_{<} R_{r, \sigma_{2}}^{-1} \right\}.  \endmultline $$
\endproclaim

\noindent
In Section 4 we are going to show that both traces in Lemma 2.5 are equal to zero by using the scattering theory.

\vskip 0.3 true cm

On the other hand, we note from (2.3) that
$$ \split
Det R_{r, \sigma_{1}}
& = Det \left( R_{r, \sigma_{2}} + \left( R_{r, \sigma_{1}} - R_{r, \sigma_{2}} \right) \right) \\
& = Det \left( R_{r, \sigma_{2}} + \frac{1}{2r} (\sigma_{2}-\sigma_{1})|_{ker B} \right) \\
& = Det \left( R_{r, \sigma_{2}} \left( I + \frac{1}{2r} R_{r, \sigma_{2}}^{-1}(\sigma_{2} - \sigma_{1})|_{ker B} \right) \right).
\endsplit
$$

\noindent
We now introduce the Fredholm determinant of a trace class operator.
Suppose that $H$ is a separable Hilbert space and $T : H \rightarrow H$ is a trace class operator. Then we define the Fredholm
determinant $det_{Fr}(I + T)$ by
$$det_{Fr}(I + T) = e^{tr \log(I + T)}.$$
The following theorem shows the relation between the Fredholm determinant and zeta-determinant ({\it cf.} Lemma 2.1 in [KV]).

\proclaim{Theorem 2.6}
$$
 \frac{Det R_{r, \sigma_{1}}}{Det R_{r, \sigma_{2}}} =
det_{Fr}\left( I + proj_{ker B} \circ \frac{1}{2r}R_{r, \sigma_{2}}^{-1}(\sigma_{2}-\sigma_{1}) \circ proj_{ker B} \right).
$$
\endproclaim
{\it Proof} : \hskip 0.3 true cm
From Lemma 2.3 and the Green theorem ({\it cf.} Lemma 4.3 in [L3]), one can see that each
$R_{r, \sigma_{i}}$ is a positive self-adjoint operator. Since every compact operator has a pure point spectrum ({\it cf.} [C]),
$ \frac{1}{2r}R_{r, \sigma_{2}}^{-1}(\sigma_{2} - \sigma_{1})|_{ker B}$ and hence
$\left( I + \frac{1}{2r}R_{r, \sigma_{2}}^{-1}(\sigma_{2} - \sigma_{1})|_{ker B} \right)$ have pure point spectra.
Since a product of positive self-adjoint operators is a positive operator,
these two facts imply that
$\left( I + \frac{1}{2r}R_{r, \sigma_{2}}^{-1}(\sigma_{2} - \sigma_{1})|_{ker B} \right)$ has only positive eigenvalues.
Hence we can choose $\pi$ as a branch-cut for logarithm.

Choose a contour $\Gamma$ in ${\Bbb C}- \{ re^{i\pi} \mid 0 \leq r < \infty \}$ containing all the eigenvalues of
$ I + \frac{1}{2r}R_{r, \sigma_{2}}^{-1}(\sigma_{2}-\sigma_{1})|_{ker B}$.
We define
$$
C := \log_{\pi} \left( I + \frac{1}{2r}R_{r, \sigma_{2}}^{-1}(\sigma_{2} - \sigma_{1})|_{ker B} \right)
\qquad \qquad \quad  \qquad \qquad  $$
$$ \qquad  = \frac{1}{2 \pi i}\int_{\Gamma} \log_{\pi} \lambda
\left( \lambda - \left( I + \frac{1}{2r}R_{r, \sigma_{2}}^{-1}(\sigma_{2}-\sigma_{1})|_{ker B} \right) \right)^{-1} d\lambda.$$
Then for $0 \leq t \leq 1$,
$$e^{tC} = \left( I + \frac{1}{2r}R_{r, \sigma_{2}}^{-1}(\sigma_{2}-\sigma_{1})|_{ker B} \right)^{t} \qquad \qquad \qquad \qquad \qquad $$
$$\quad = \frac{1}{2\pi i} \int_{\Gamma} \lambda^{t} \left( \lambda - \left( I + \frac{1}{2r}R_{r, \sigma_{2}}^{-1}
(\sigma_{2}-\sigma_{1})|_{ker B} \right) \right)^{-1} d\lambda.$$
We put
$$R_{r,\sigma_{2}}(t) = R_{r,\sigma_{2}}
\left( I + \frac{1}{2r}R_{r, \sigma_{2}}^{-1} (\sigma_{2} - \sigma_{1})|_{ker B} \right)^{t} = R_{r,\sigma_{2}} e^{tC} .$$
Since $\left( I + \frac{1}{2r}R_{r, \sigma_{2}}^{-1}(\sigma_{2} - \sigma_{1})|_{ker B} \right)^{t}$ has only positive eigenvalues,
$R_{r,\sigma_{2}}(t)$ also has only positive eigenvalues. Hence we can also choose $\pi$ as a branch-cut for logarithm concerning
the operator $R_{r,\sigma_{2}}(t)$.
For $Re s \gg 0$ we have
$$
\zeta_{R_{r, \sigma_{2}}(t)}(s) := Tr \left( \frac{1}{2\pi i}
\int_{\gamma}\lambda^{-s} (\lambda - R_{r,\sigma_{2}}(t))^{-1} d\lambda \right),$$
where for small $\epsilon > 0$
$$\gamma = \{ re^{i\frac{\pi}{4}} \mid \infty > r \geq \epsilon \} \cup
\{ \epsilon e^{i\phi} \mid \frac{\pi}{4} \geq \phi \geq -\frac{\pi}{4} \} \cup
\{ re^{-i\frac{\pi}{4}} \mid \epsilon \leq r < \infty \}.$$
Then for $Re s \gg 0$,
$$\split
\frac{d}{dt} \zeta_{R_{r,\sigma_{2}}(t)}(s)
& = Tr \left( \frac{1}{2\pi i} \int_{\gamma} \lambda^{-s} (\lambda - R_{r,\sigma_{2}}(t))^{-1}  R_{r,\sigma_{2}}(t) C
(\lambda - R_{r, \sigma_{2}}(t))^{-1} d\lambda \right) \\
& = Tr \left( R_{r,\sigma_{2}}(t) C \frac{1}{2\pi i} \int_{\gamma} \lambda^{-s} (\lambda - R_{r,\sigma_{2}}(t))^{-2}
 d\lambda \right) \\
& = Tr \left( R_{r,\sigma_{2}}(t) C \frac{-1}{2\pi i} \int_{\gamma} \lambda^{-s}
\frac{\lambda}{d\lambda} (\lambda - R_{r,\sigma_{2}}(t))^{-1}  d\lambda \right) \\
& = -s Tr \left( R_{r,\sigma_{2}}(t) C \frac{1}{2\pi i} \int_{\gamma} \lambda^{-s-1}
(\lambda - R_{r,\sigma_{2}}(t))^{-1}  d\lambda \right) \\
& = - s Tr \left( R_{r,\sigma_{2}}(t) C ( R_{r,\sigma_{2}}(t) )^{-s-1} \right) \\
& = - s Tr \left( C  ( R_{r,\sigma_{2}}(t) )^{-s} \right).
\endsplit
$$
Taking derivative with respect to $s$ at $s=0$ gives
$$\frac{d}{ds}\frac{d}{dt} \zeta_{R_{r,\sigma_{2}}(t)}(s) |_{s=0} =
\frac{d}{dt}\frac{d}{ds} \zeta_{R_{r,\sigma_{2}}(t)}(s) |_{s=0} = - Tr (C) .$$
Integrating from 0 to 1, we have
$$
\log Det R_{r, \sigma_{2}} (1) - \log Det R_{r, \sigma_{2}} (0) =  Tr (C) $$
and therefore
$$ \split
\frac{Det R_{r, \sigma_{1}}}{Det R_{r, \sigma_{2}}}
& = e^{Tr (C)} = e^{Tr \log_{\pi}\left( I + \frac{1}{2r}R_{r, \sigma_{2}}^{-1}(\sigma_{2}-\sigma_{1})|_{ker B} \right)} \\
& = det_{Fr} e^{\log_{\pi}\left( I + \frac{1}{2r}R_{r, \sigma_{2}}^{-1}(\sigma_{2} - \sigma_{1})|_{ker B} \right)} \\
& = det_{Fr} \left( I + \frac{1}{2r}R_{r, \sigma_{2}}^{-1}(\sigma_{2} - \sigma_{1})|_{ker B} \right)  \qquad \qquad \qquad \\
& = det_{Fr} \left( I +  proj_{ker B} \circ \frac{1}{2r}R_{r, \sigma_{2}}^{-1}(\sigma_{2} - \sigma_{1}) \circ proj_{ker B} \right).
\qquad \qquad  \qed
\endsplit
$$
In Section 4 we are going to show that $\frac{Det R_{r, \sigma_{1}}}{Det R_{r, \sigma_{2}}}$ does not
depend on the cylinder length $r$
and obtain Theorem 1.3 by computing the limit
$$\lim_{r\rightarrow\infty}
det_{Fr} \left( I + proj_{ker B} \circ \frac{1}{2r}R_{r, \sigma_{2}}^{-1}(\sigma_{2} - \sigma_{1}) \circ proj_{ker B} \right) .
\tag2.5 $$

\vskip 1 true cm

\S 3 {\bf W. M\"uller's scattering theory }

\vskip 0.5 true cm

In this section we review some basic facts of the scattering theory developed by W. M\"uller in [M] which is needed in this paper.
Recall that $M_{\infty} = M \cup_{Y} [0, \infty)$ and $D_{M_{\infty}}$, the natural extension of $D_{M}$ to $M_{\infty}$,
has the form $G(\partial_{u}+B)$ on the cylinder part.
Let $\mu_{1}$ be the smallest positive eigenvalue of $B$. Then it is shown in [M] that
for $\lambda \in {\Bbb R}$ with $|\lambda| < \mu_{1}$
there exists a regular one-parameter family of unitary operators $C(\lambda)$ on $ker B$, called the scattering matrices, satisfying
the following properties :
\newline
$$ (1) \quad C(\lambda) C(-\lambda) = I  \quad \text{ and } \quad  C(\lambda) G = - G C(\lambda)  \qquad \qquad
\qquad \qquad \qquad \qquad \qquad  \tag3.1 $$
\noindent
Putting $\lambda = 0$, we have
$$ C(0)^{2} = I \quad  \text{ and } \quad C(0) G = - G C(0).$$
Hence, $C(0)$ is the natural choice of a unitary involution anticommuting with $G$ on $ker B$.
\newline
(2) For any $f \in ker B$ and $\lambda \in {\Bbb R}$ with $|\lambda| < \mu_{1}$,
the generalized eigensection $E(f,\lambda)$ of
$D_{M_{\infty}}$ attached to $f$ can be expressed, on the cylinder part, by
$$E(f,\lambda) = e^{- i \lambda u} (f - i G f) + e^{i \lambda u} C(\lambda) (f - i G f) + \theta(f, \lambda),
\tag3.2 $$
where $\theta(f, \lambda)$ is square integrable and $\theta(f, \lambda, (u, \cdot ))$ is orthogonal to $ker B$.
Moreover, $E(f,\lambda) $ satisfies
$D_{M_{\infty}} E(f,\lambda) = \lambda E(f,\lambda) $.

\vskip 0.3 true cm

Now we suppose that $f \in ker (I - C(0))$. Then we have

$$ \frac{1}{2} E(f, 0) = f + \frac{1}{2} \theta(f,0). \tag3.3 $$
\proclaim{Definition 3.1}
We call $E(f, 0)$ an extended $L^{2}$-solution of  $D_{M_{\infty}}$ and $f$ the limiting value of $ \frac{1}{2} E(f, 0)$.
Hence the dimension of the space of limiting values is $\frac{1}{2} dim ker B$.
\endproclaim
Let $f$ belong to $ker (I - C(0))$ and $E(f,\lambda)$ be the generalized eigensection attached to $f$ as (3.2).
Then we can compute the $L^{2}$-part $\theta(f, \lambda)$
on the cylinder part as follows.
Put
$$ \theta(f, \lambda) = \sum_{0 < \mu_{j} \in Spec(B)} (p_{j}(u) \phi_{\mu_{j}} + q_{j}(u) \phi_{-\mu_{j}}),$$
where $\phi_{\mu_{j}}$ ($\phi_{-\mu_{j}} = G\phi_{\mu_{j}}$) is an eigensection of $B$
corresponding to the eigenvalue $\mu_{j}$ ($-\mu_{j}$).
Since $G(\partial_{u}+B) \theta(f, \lambda) = \lambda \theta(f, \lambda)$ and $\theta(f, \lambda)$ satisfies $L^{2}$-condition,
direct computations show that
$$ \left\{\aligned
p_{j}(u) & = a_{j}(\lambda) ( \sqrt{\mu_{j}^{2} + \lambda^{2}} + \mu_{j}) e^{-\sqrt{\mu_{j}^{2} + \lambda^{2}} u} \\
q_{j}(u) & = - \lambda a_{j}(\lambda) e^{-\sqrt{\mu_{j}^{2} + \lambda^{2}} u},
\endaligned \right.
$$
where $a_{j}(\lambda)$ is a smooth function of $\lambda$ ($|\lambda| < \mu_{1}$). Hence we have
$$
E(f,\lambda) =  e^{- i \lambda u} (f - i G f) + e^{i \lambda u} C(\lambda) (f - i G f) + $$
$$
\sum_{0 < \mu_{j} \in Spec(B)} \left( a_{j}(\lambda) ( \sqrt{\mu_{j}^{2} + \lambda^{2}} + \mu_{j}) e^{-\sqrt{\mu_{j}^{2} + \lambda^{2}} u}
\phi_{\mu_{j}} - \lambda a_{j}(\lambda) e^{-\sqrt{\mu_{j}^{2} + \lambda^{2}} u}\phi_{-\mu_{j}} \right) \tag3.4
$$
and
$$E(f,0) = 2f + \sum_{0 < \mu_{j} \in Spec(B)} 2\mu_{j} a_{j}(0) e^{-\mu_{j} u} \phi_{\mu_{j}}. \tag3.5 $$
Since $C(0)f = f$ and $C^{\prime}(0) G = - G C^{\prime}(0)$,
the derivative of (3.4) with respect to $\lambda$ at $\lambda = 0$ is,
on the cylinder part,
$$ \multline
\frac{d}{d\lambda} E(f,\lambda)|_{\lambda = 0} = -2 u G f + C^{\prime}(0) (f - i G f) + \\
 \sum_{0 < \mu_{j} \in Spec(B)} \left( 2\mu_{j}a_{j}^{\prime}(0) e^{-\mu_{j} u}\phi_{\mu_{j}} -
a_{j}(0) e^{-\mu_{j} u} \phi_{-\mu_{j}} \right). \endmultline   \tag3.6
$$
Note that even if $\frac{d}{d\lambda} E(f,\lambda)|_{\lambda = 0}$ is not a solution of $D_{M_{\infty}}$,
it is a solution of $D^{2}_{M_{\infty}}$. Hence,
if we denote $l = dim ker (I - C(0))$, the equation (3.5) and (3.6) provide $2l$ distinct solutions
of $D^{2}_{M_{\infty}}$ on $M_{\infty}$. Using this observation we are going to prove Theorem 1.2 and 1.3
 in the next section.

\vskip 1 true cm

\S 4 {\bf The proof of Theorem 1.2  and 1.3 }

\vskip 0.5 true cm

In this paper we assume that both $D_{M, P_{\sigma_{1}}}^{2}$ and $D_{M, P_{\sigma_{2}}}^{2}$ are invertible operators.
This assumption implies that both $D_{M_{r},P_{\sigma_{1}}}^{2}$ and $D_{M_{r},P_{\sigma_{2}}}^{2}$ are invertible for each $r > 0$
and that
$$
ker (I - \sigma_{1}) \cap ker (I - C(0))  = \{ 0 \} = ker (I - \sigma_{2}) \cap ker (I - C(0)) .   \tag4.1 $$

Suppose that $\{ f_{1}, f_{2}, \cdots, f_{l} \}$ is an orthonormal basis for $ker (I - C(0))$.
Then for each $k$ ($1 \leq k \leq l$) the extended $L^{2}$- solution $\frac{1}{2}E(f_{k},0)$ attached to $f_{k}$
has the form, on the cylinder part,
$$\frac{1}{2}E(f_{k},0) = f_{k} + \sum_{0<\mu_{j}\in Spec(B)}b^{k}_{j}e^{-\mu_{j} u}\phi_{\mu_{j}}.$$
For each section $f \in C^{\infty}(Y_{r})$ we denote
$f^{+} = \frac{I + \sigma_{2}}{2}(f)$ and $f^{-} = \frac{I - \sigma_{2}}{2}(f)$.
Then we can decompose $f_{k}$ into $f_{k}=f_{k}^{+} + f_{k}^{-}$, where $\sigma_{2}f_{k}^{+} = f_{k}^{+}$ and
$\sigma_{2}f_{k}^{-} = - f_{k}^{-}$.
Define
$$\psi_{k} = \cases
\frac{1}{2}E(f_{k},0) & \text{ on $M$ }  \\
f_{k} - \frac{u}{r}f_{k}^{-} + \sum_{\mu_{j}>0}b_{j}^{k}e^{-\mu_{j} u}\phi_{\mu_{j}} & \text{ if } 0 \leq u \leq r .
\endcases
$$
Then $\psi_{k}$ satisfies the boundary condition $P_{\sigma_{2}}$ on $Y_{r}$.
Hence
$$ R_{r, \sigma_{2}}\left(f_{k} + \sum_{\mu_{j}>0}b_{j}^{k}\phi_{\mu_{j}}\right)  = $$
$$\frac{\partial}{\partial u} \left(\frac{1}{2}E(f_{k},0)\right)|_{u=0} -
\frac{\partial}{\partial u} \left(f_{k} - \frac{u}{r}f_{k}^{-} + \sum_{\mu_{j}>0}b_{j}^{k}e^{-\mu_{j} u}\phi_{\mu_{j}}\right)|_{u=0}$$
$$ =  \frac{1}{r}f_{k}^{-},$$
and
$$
\frac{1}{2r}R_{r, \sigma_{2}}^{-1}(f_{k}^{-}) =
\frac{1}{2}f_{k} + \frac{1}{2}\sum_{\mu_{j}>0}b_{j}^{k}\phi_{\mu_{j}}. \tag4.2 $$

We denote by $P$ the orthogonal projection from $ker (I - \sigma_{2})$ to $ker (I + C(0))$ {\it i.e.}
$$ P = \frac{I - C(0)}{2} \frac{I + \sigma_{2}}{2} : ker (I - \sigma_{2}) \rightarrow ker (I + C(0)). $$

The equality (4.1) implies that $P$ is invertible and we denote by $P^{-1}$ the inverse of $P$.
Note that
$$ f_{k}^{-} = \frac{I + C(0)}{2} f_{k}^{-} + \frac{I - C(0)}{2} f_{k}^{-}. $$
Since
$$
P^{-1}\left(\frac{I - C(0)}{2} f_{k}^{-}\right) =
\frac{I + C(0)}{2}P^{-1}\left(\frac{I - C(0)}{2} f_{k}^{-}\right) + \frac{I - C(0)}{2} f_{k}^{-},$$
$$
f_{k}^{-} - P^{-1}\left(\frac{I - C(0)}{2} f_{k}^{-} \right) \in ker(I - C(0)). $$
The equality (4.1) also implies that the orthogonal projection $\frac{I - \sigma_{2}}{2}$
from $ker (I - C(0))$ to $ker(I + \sigma_{2})$ is an isomorphism.
This fact and the following equality
$$
\frac{I - \sigma_{2}}{2} f_{k} =
\frac{I - \sigma_{2}}{2}\left(f_{k}^{-} - P^{-1}\left(\frac{I - C(0)}{2} f_{k}^{-}\right)\right) = f_{k}^{-} $$
imply that
$$ f_{k} = f_{k}^{-} - P^{-1}\left(\frac{I - C(0)}{2} f_{k}^{-}\right). \tag4.3  $$
The equation (4.2) and (4.3) lead to
$$
\frac{1}{2r} R^{-1}_{r,\sigma_{2}}(f_{k}^{-}) = \frac{1}{2}\left(I
- P^{-1}\left(\frac{I - C(0)}{2}\right)\right) f_{k}^{-} +
\frac{1}{2} \sum_{\mu_{j}>0} b_{j}^{k}\phi_{\mu_{j}}.  \tag4.4 $$
Summarizing the above argument, we have the following lemma.
\proclaim{Lemma 4.1}
$$
proj_{ker B} \circ \frac{1}{2r}R^{-1}_{r,\sigma_{2}} \circ proj_{ker(I + \sigma_{2})} =
\frac{1}{2}\left(I - P^{-1}\left(\frac{I - C(0)}{2}\right)\right) \frac{I - \sigma_{2}}{2}.
$$
\endproclaim

Before describing $proj_{ker B} \circ \frac{1}{2r}R^{-1}_{r,\sigma_{2}} \circ proj_{ker(I - \sigma_{2})}$
we prove Theorem 1.2 by using the equation (4.2), (4.4) and Lemma 2.5.
We recall :

\proclaim{Theorem 1.2}
Suppose that both $D_{M, P_{\sigma_{1}}}^{2}$ and $D_{M, P_{\sigma_{2}}}^{2}$ are invertible operators. Then
$$\frac{Det D_{M_{r}, P_{\sigma_{1}}}^{2}}{Det D_{M_{r}, P_{\sigma_{2}}}^{2}}$$
does not depend on the cylinder length $r$.
\endproclaim
{\it Proof} : \hskip 0.3 true cm
By Lemma 2.5 it's enough to show that
$$ Tr \left\{ proj_{ker(I + \sigma_{2})} R_{r, \sigma_{2}}^{-1} proj_{ker(I + \sigma_{2})} -
proj_{ker(I + \sigma_{1})} R_{r, \sigma_{1}}^{-1} proj_{ker(I + \sigma_{1})} \right\} = 0, $$
$$ Tr \left\{ (\sigma_{2} - \sigma_{1})|_{ker B} R_{r, \sigma_{1}}^{-1}
K_{r} P_{<} R_{r, \sigma_{2}}^{-1} \right\} = 0.$$
\noindent
The equation (4.2) shows that
$proj_{ker(I + \sigma_{i})} R_{r, \sigma_{i}}^{-1} proj_{ker(I + \sigma_{i})}$ ($i = 1, 2$)
has only one eigenvalue $r$ of multiplicity
$\frac{1}{2} dim ker B$, which implies the first equality.

For the second equality we suppose that $\{ g_{1}, g_{2}, \cdots, g_{l} \}$  is an orthonormal basis for $ker(I + \sigma_{2})$.
Then $\{ g_{1}, g_{2}, \cdots, g_{l}, Gg_{1}, Gg_{2}, \cdots, Gg_{l} \}$ is an orthonormal basis for $ker B$.
Hence,
$$ \multline
Tr \left\{ (\sigma_{2} - \sigma_{1})|_{ker B} R_{r, \sigma_{1}}^{-1}
K_{r} P_{<} R_{r, \sigma_{2}}^{-1} \right\}
 = \sum_{i=1}^{l} \langle (\sigma_{2} - \sigma_{1})|_{ker B} R_{r, \sigma_{1}}^{-1}
K_{r} P_{<} R_{r, \sigma_{2}}^{-1} g_{i}, g_{i} \rangle \\
 + \sum_{i=1}^{l} \langle (\sigma_{2} - \sigma_{1})|_{ker B} R_{r, \sigma_{1}}^{-1}
K_{r} P_{<} R_{r, \sigma_{2}}^{-1} Gg_{i}, Gg_{i} \rangle . \endmultline $$

\noindent
The equation (4.4) shows that $P_{<} R_{r, \sigma_{2}}^{-1} g_{i} = 0$.
Since $(\sigma_{2} - \sigma_{1})|_{ker B}$, $R_{r, \sigma_{1}}^{-1}$ and
$K_{r} P_{<}$ are self-adjoint operators, we have
$$ \split
\langle (\sigma_{2} - \sigma_{1})|_{ker B} R_{r, \sigma_{1}}^{-1}
K_{r} P_{<} R_{r, \sigma_{2}}^{-1} Gg_{i}, Gg_{i} \rangle
& = \langle R_{r, \sigma_{2}}^{-1} Gg_{i}, K_{r} P_{<} R_{r, \sigma_{1}}^{-1} (\sigma_{2} - \sigma_{1}) Gg_{i} \rangle \\
& = \langle R_{r, \sigma_{2}}^{-1} Gg_{i}, K_{r} P_{<} R_{r, \sigma_{1}}^{-1} ( I - \sigma_{1}) Gg_{i} \rangle .
\endsplit $$

\noindent
Since $( I - \sigma_{1}) Gg_{i}$ belongs to $ker(I + \sigma_{1})$, the equation (4.4) again shows that
$P_{<} R_{r, \sigma_{1}}^{-1} ( I - \sigma_{1}) Gg_{i} = 0$.
This completes the proof of the theorem.
\qed

\vskip 0.5 true cm

Now we describe the operator $proj_{ker B} \circ \frac{1}{2r}R^{-1}_{r,\sigma_{2}} \circ proj_{ker(I - \sigma_{2})}$.
We begin with the equation (3.6), which is at $u=0$
$$
\frac{d}{d\lambda} E(f_{k},\lambda)|_{\lambda = 0, u=0} = C^{\prime}(0) (f_{k} - i G f_{k}) +
\sum_{0 < \mu_{j} \in Spec(B)} \left( 2\mu_{j}a_{j}^{\prime}(0) \phi_{\mu_{j}} -
a_{j}(0) \phi_{-\mu_{j}} \right).   \tag4.5
$$
To compute $R_{r,\sigma_{2}}^{-1}(\frac{d}{d\lambda} E(f_{k},\lambda)|_{\lambda = 0, u=0})$
we define a section $\psi_{k}$ as follows.
First, on $[0, r]\times Y$ we put
$$
\psi_{k,cyl} = C^{\prime}(0)(f_{k} - i G f_{k}) - \frac{u}{r} \alpha_{k}^{-} +
\sum_{\mu_{j}>0} 2\mu_{j}a_{j}^{\prime}(0)e^{-\mu_{j} u}\phi_{\mu_{j}} $$
$$ - \sum_{\mu_{j}>0} a_{j}(0)\frac{e^{-\mu_{j}(u-r)} - e^{\mu_{j}(u-r)}}{e^{\mu_{j}r} - e^{-\mu_{j}r}} \phi_{-\mu_{j}},
$$
where $\alpha_{k}^{-} = \frac{I - \sigma_{2}}{2}(C^{\prime}(0)(f_{k} - i G f_{k}))$.
Now we define $\psi_{k}$ by
$$
\psi_{k} = \cases
\frac{d}{d\lambda} E(f_{k},\lambda)|_{\lambda=0}  & \text{ on } M  \\
\psi_{k,cyl}  & \text{ on } [0, r] \times Y .
\endcases
$$
Then $\psi_{k}$ is a continuous section with $\psi_{k}|_{Y_{0}} = \frac{d}{d\lambda} E(f_{k},\lambda)|_{\lambda=0, u=0}$
and satisfies the boundary condition $P_{\sigma_{2}}$ on $Y_{r}$.
Furthermore, $D^{2}_{M_{r}}\psi_{k}=0$ on $M_{r}-Y_{0}$.
Hence
$$
R_{r,\sigma_{2}}\left(\frac{d}{d\lambda} E(f_{k},\lambda)|_{\lambda=0, u=0}\right) =
\frac{\partial}{\partial u}\left(\frac{d}{d\lambda} E(f_{k},\lambda)|_{\lambda=0}\right)|_{u=0} -
\frac{\partial}{\partial u}(\psi_{k,cyl})|_{u=0} $$
$$
= -2 G f_{k} + \frac{1}{r} \alpha_{k}^{-} -
\sum_{\mu_{j}>0}\frac{2\mu_{j}e^{-\mu_{j}r}}{e^{\mu_{j}r} - e^{-\mu_{j}r}}a_{j}(0)\phi_{-\mu_{j}}.  \tag4.6 $$
The equation (4.5) and (4.6) lead to
$$
R_{r,\sigma_{2}}^{-1}(-2 G f_{k} + \frac{1}{r}\alpha_{k}^{-}) = C^{\prime}(0)(f_{k} - i G f_{k}) +
\sum_{\mu_{j}>0}(2\mu_{j} a_{j}^{\prime}(0)\phi_{\mu_{j}} - a_{j}(0)\phi_{-\mu_{j}}) $$
$$
+ R_{r,\sigma_{2}}^{-1}\left(\sum_{\mu_{j}>0}\frac{2\mu_{j}e^{-\mu_{j}r}}{e^{\mu_{j}r}-e^{-\mu_{j}r}}a_{j}(0)\phi_{-\mu_{j}}\right).
\tag4.7 $$

\noindent
Since $C^{\prime}(0) f_{k} \in ker(I-C(0))$,   the equation (4.2) shows that
$$
\frac{1}{r}R_{r,\sigma_{2}}^{-1}((C^{\prime}(0) f_{k})^{-}) =
C^{\prime}(0) f_{k} + \sum_{\mu_{j}>0}{\tilde b_{j}}\phi_{\mu_{j}}. \tag4.8 $$
Combining (4.7), (4.8) with Lemma 4.1, we have
$$
proj_{ker B} \circ R_{r,\sigma_{2}}^{-1}(- 2 G f_{k} )  $$
$$ \multline
 = - proj_{ker B} \circ \frac{1}{r} R_{r,\sigma_{2}}^{-1} ((i G C^{\prime}(0) f_{k})^{-} ) +i G C^{\prime}(0) f_{k}  \\
+ proj_{ker B} \circ R_{r,\sigma_{2}}^{-1}
\left(\sum_{\mu_{j}>0} \frac{2\mu_{j}e^{-\mu_{j}r}}{e^{\mu_{j}r}-e^{-\mu_{j}r}} a_{j}(0)\phi_{-\mu_{j}}\right)
\endmultline
$$
$$ \multline
= - \left(I - P^{-1}\left(\frac{I - C(0)}{2}\right)\right)\left(\frac{I - \sigma_{2}}{2}\right)(i G C^{\prime}(0) f_{k}) \\
+  i G C^{\prime}(0) f_{k} + proj_{ker B} \circ R_{r,\sigma_{2}}^{-1}
\left(\sum_{\mu_{j}>0} \frac{2\mu_{j}e^{-\mu_{j}r}}{e^{\mu_{j}r}-e^{-\mu_{j}r}} a_{j}(0)\phi_{-\mu_{j}}\right).
\endmultline \tag4.9
$$
Using (4.9) and Lemma 4.1 again, we have
$$
proj_{ker B} \circ \frac{1}{2r} R_{r,\sigma_{2}}^{-1}(G f_{k})^{+}  =  proj_{ker B} \circ \frac{1}{2r}R_{r,\sigma_{2}}^{-1}(G f_{k}) -
proj_{ker B} \circ \frac{1}{2r}R_{r,\sigma_{2}}^{-1}(G f_{k})^{-} $$
$$ \multline
= \frac{1}{4r} \left(I - P^{-1}\left(\frac{I - C(0)}{2}\right)\right)\left(\frac{I - \sigma_{2}}{2}\right)(i G C^{\prime}(0) f_{k}) \\
- \frac{i}{4r} G C^{\prime}(0) f_{k}  - \frac{1}{4r} proj_{ker B} \circ R_{r,\sigma_{2}}^{-1}
\left(\sum_{\mu_{j}>0} \frac{2\mu_{j}e^{-\mu_{j}r}}{e^{\mu_{j}r}-e^{-\mu_{j}r}} a_{j}(0)\phi_{-\mu_{j}}\right)
\endmultline
$$
$$ - \frac{1}{2} \left( I - P^{-1}\left(\frac{I - C(0)}{2}\right)\right) \left(\frac{I - \sigma_{2}}{2}\right)(G f_{k}).$$
Therefore, we have
$$
\lim_{r\rightarrow\infty} proj_{ker B} \circ \frac{1}{2r} R_{r,\sigma_{2}}^{-1}(G f_{k})^{+} = -\frac{1}{2}
\left( I - P^{-1}\left(\frac{I - C(0)}{2}\right)\right) \left(\frac{I - \sigma_{2}}{2}\right)(G f_{k}). \tag4.10 $$

\vskip 0.3 true cm

Now we denote by $K$ the orthogonal projection from $ker(I + \sigma_{2})$ onto $ker(I - C(0))$,
{\it i.e.}
$$K=\frac{I + C(0)}{2} \frac{I - \sigma_{2}}{2}. $$
The equality (4.1) implies that $K$ is invertible and we denote the inverse by $K^{-1}$.
We decompose $(G f_{k})^{+}$ by
$$
(G f_{k})^{+} = \frac{I + C(0)}{2} (G f_{k})^{+} + \frac{I - C(0)}{2} (G f_{k})^{+}.$$
Since
$$
K^{-1}\left( \frac{I + C(0)}{2} (G f_{k})^{+} \right) \qquad \qquad \qquad $$
$$ = \frac{I + C(0)}{2} K^{-1}\left( \frac{I + C(0)}{2} (G f_{k})^{+} \right) +
\frac{I - C(0)}{2} K^{-1}\left( \frac{I + C(0)}{2} (G f_{k})^{+} \right) $$
$$
= \frac{I + C(0)}{2} (G f_{k})^{+} + \frac{I - C(0)}{2} K^{-1}\left( \frac{I + C(0)}{2} (G f_{k})^{+} \right) ,
\qquad \qquad \qquad \quad $$
we have
$$
(G f_{k})^{+} - K^{-1}\left( \frac{I + C(0)}{2} (G f_{k})^{+} \right) \in ker (I + C(0)) $$
and
$$
\frac{I + \sigma_{2}}{2} \left((G f_{k})^{+} - K^{-1}\left( \frac{I + C(0)}{2} (G f_{k})^{+} \right)\right)
= (G f_{k})^{+}.$$
Since  $\frac{I + \sigma_{2}}{2} : ker (I + C(0)) \rightarrow ker( I - \sigma_{2}) $ is an isomorphism (see (4.1))
and $\frac{I + \sigma_{2}}{2} (G f_{k}) = (G f_{k})^{+}$,
$$\split
G f_{k} & = (G f_{k})^{+} - K^{-1}\left( \frac{I + C(0)}{2} (G f_{k})^{+} \right) \\
& = \left(I - K^{-1}\left( \frac{I + C(0)}{2} \right)\right) (G f_{k})^{+} .
\endsplit   \tag4.11
$$
Therefore, the equation (4.10) and (4.11) lead to the following result.
\proclaim{Lemma 4.2}
$$
\lim_{r\rightarrow\infty} proj_{ker B}\circ \frac{1}{2r}R^{-1}_{\sigma_{2}} \circ proj_{ker(I-\sigma_{2})} $$
$$= -\frac{1}{2}\left( I- P^{-1}\left( \frac{I-C(0)}{2}\right)\right) \frac{I-\sigma_{2}}{2}
\left( I- K^{-1}\left( \frac{I+C(0)}{2}\right)\right) \frac{I+\sigma_{2}}{2} $$
$$ = \frac{1}{2}\left( I- P^{-1}\left( \frac{I-C(0)}{2}\right)\right) \frac{I-\sigma_{2}}{2}
K^{-1}\left( \frac{I+C(0)}{2}\right) \frac{I+\sigma_{2}}{2}. \qquad \qquad $$
\endproclaim

\vskip 0.3 true cm

We combine Lemma 4.1 and 4.2 to have the following.
$$
\lim_{r\rightarrow\infty} proj_{ker B}\circ \frac{1}{2r}R^{-1}_{\sigma_{2}} \circ proj_{ker B} $$
$$ \multline
= \frac{1}{2}\left( I- P^{-1}\left( \frac{I-C(0)}{2}\right)\right) \frac{I-\sigma_{2}}{2} \\
+ \frac{1}{2}\left( I- P^{-1}\left( \frac{I-C(0)}{2}\right)\right) \frac{I-\sigma_{2}}{2}
K^{-1}\left( \frac{I+C(0)}{2}\right) \frac{I+\sigma_{2}}{2}
\endmultline $$
$$
= \frac{1}{2}\left( I- P^{-1}\left( \frac{I-C(0)}{2}\right)\right) \frac{I-\sigma_{2}}{2}
\left\{ I + K^{-1}\left( \frac{I+C(0)}{2}\right) \frac{I+\sigma_{2}}{2} \right\}.  \qquad $$

\vskip 0.3 true cm

Now we define $T : ker (I - \sigma_{2}) \oplus  ker (I + \sigma_{2}) \rightarrow ker ( I + C(0) ) \oplus ker ( I - C(0) ) $ by
$$ T = \frac{I - C(0)}{2} \frac{I+\sigma_{2}}{2} +  \frac{I + C(0)}{2} \frac{I-\sigma_{2}}{2}.$$
In other words, if $f^{+} \in ker(I - \sigma_{2})$, $T(f^{+}) = P(f^{+})$ and
if $f^{-} \in ker(I + \sigma_{2})$, $T(f^{-}) = K(f^{-})$ .
Note that $T$ can be rewritten as
$$
T = \frac{1}{2} (I - C(0)\sigma_{2}) = \frac{1}{2} C(0) (C(0) - \sigma_{2}) = \frac{1}{2} (\sigma_{2} - C(0)) \sigma_{2}.
\tag4.12 $$
We use (4.12) to obtain the following equalities.
$$
\lim_{r\rightarrow\infty} proj_{ker B}\circ \frac{1}{2r}R^{-1}_{\sigma_{2}} \circ proj_{ker B} $$
$$
= \frac{1}{2}\left( I- T^{-1}\left( \frac{I-C(0)}{2}\right)\right) \frac{I - \sigma_{2}}{2}
\left\{ \frac{I-\sigma_{2}}{2} + T^{-1}\left( \frac{I+C(0)}{2}\right) \frac{I+\sigma_{2}}{2} \right\} \qquad \qquad $$
$$
= \frac{1}{2} T^{-1} \left( T -  \frac{I-C(0)}{2}\right) \frac{I - \sigma_{2}}{2}
T^{-1} \left\{ T \frac{I-\sigma_{2}}{2} +  \frac{I+C(0)}{2} \frac{I+\sigma_{2}}{2} \right\} \qquad\qquad\qquad $$
$$
= \frac{1}{2} T^{-1} \left( \frac{I + C(0)}{2} -  \frac{I - C(0)}{2}\right) \frac{I - \sigma_{2}}{2}
T^{-1} \left\{ \frac{I + C(0)}{2} \frac{I - \sigma_{2}}{2} +  \frac{I + C(0)}{2}
 \frac{I + \sigma_{2}}{2} \right\} $$
$$
= \frac{1}{2} T^{-1} C(0) \frac{I - \sigma_{2}}{2} T^{-1} \frac{I + C(0)}{2} \qquad\qquad\qquad\qquad\qquad
\qquad\qquad\qquad\qquad\qquad \quad $$
$$
= \frac{1}{2} T^{-1} C(0)  T^{-1} \frac{I + C(0)}{2} .
\qquad\qquad\qquad\qquad\qquad\qquad\qquad\qquad\qquad\qquad \qquad \qquad  \tag4.13 $$

\noindent
The equation (4.12) shows that
$$ T^{-1} = 2 (C(0) - \sigma_{2})^{-1} C(0)$$
and hence
$$ T^{-1} C(0) T^{-1} = 4 (C(0) - \sigma_{2})^{-2} C(0). \qquad \qquad  \tag4.14 $$
Combining (4.13) and (4.14) we have
$$
\lim_{r\rightarrow\infty} det
\left( I + proj_{ker B} \circ \frac{1}{2r} R^{-1}_{\sigma_{2}}  (\sigma_{2} - \sigma_{1}) \circ proj_{kerB} \right)  $$
$$ = det \left( I +   (C(0) - \sigma_{2})^{-2} (I + C(0))  (\sigma_{2} - \sigma_{1}) \right)  $$
$$ \qquad\quad\quad = det \left( I +  (I + C(0)) (\sigma_{2} - \sigma_{1})  (C(0) - \sigma_{2})^{-2} (I + C(0)) \right)  $$
$$ = det \left( I +  (\sigma_{2} - \sigma_{1})  (C(0) - \sigma_{2})^{-2} (I + C(0)) \right) .  \tag4.15 $$

\vskip 0.3 true cm

In the next step we are going to make the formula (4.15) as simple as possible.
The following lemma and corollary are easy to prove but very useful.
\proclaim{Lemma 4.3}
For each $i = 1, 2$,
$$
C(0) \left( C(0) - \sigma_{i} \right)^{-1} = - \left( C(0) - \sigma_{i} \right)^{-1} \sigma_{i}  \quad \text{ and } $$
$$ \left( C(0) - \sigma_{i} \right)^{-1} C(0) = - \sigma_{i} \left( C(0) - \sigma_{i} \right)^{-1} . \qquad \quad $$
\endproclaim
{\it Proof} : \hskip 0.3 true cm
The lemma follows from the relations
$$ \left( C(0) - \sigma_{i} \right) C(0) = \sigma_{i} \left( \sigma_{i} - C(0) \right) \text{ and }
 C(0) \left( C(0) - \sigma_{i} \right) = \left( \sigma_{i} - C(0) \right) \sigma_{i}. $$
\qed

Using Lemma 4.3 twice leads to the following corollary.

\proclaim{Corollary 4.4}
For each $i = 1, 2$,
$$
C(0) \left( C(0) - \sigma_{i} \right)^{-2} =  \left( C(0) - \sigma_{i} \right)^{-2} C(0)  \quad \text{ and } $$
$$ \left( C(0) - \sigma_{i} \right)^{-2} \sigma_{i} = \sigma_{i} \left( C(0) - \sigma_{i} \right)^{-2} . \qquad \quad $$
\endproclaim

Now for $0 < t < 1$ we put
$$A_{\sigma_{i}}(t) = \frac{1}{2} \left( I - C(0) \right) + \frac{t}{2} \left( I - \sigma_{i} \right).$$
One can check easily by using Corollary 4.4 that
$$
A_{\sigma_{i}}(t)^{-1} = \frac{2}{t} \left( C(0) - \sigma_{i} \right)^{-2} \left( I + C(0) \right)
+ 2 \left( C(0) - \sigma_{i} \right)^{-2} \left( I + \sigma_{i} \right). $$
Using Corollary 4.4 again we have
$$ \multline
A_{\sigma_{1}}(t) A_{\sigma_{2}}(t)^{-1} =
I + \left( \sigma_{2} - \sigma_{1}\right)  \left( C(0) - \sigma_{2}\right)^{-2} \left( C(0) + I \right) \\
+ t \left( I - \sigma_{1} \right) \left( C(0) - \sigma_{2} \right)^{-2} \left( I + \sigma_{2} \right) .
\endmultline \tag4.16 $$

Therefore, Theorem 1.2, Theorem 2.6 and the equation (4.15), (4.16) imply the following equalities.
$$ \split
\frac{Det D^{2}_{M, P_{\sigma_{1}}}}{Det D^{2}_{M, P_{\sigma_{2}}}} & =
\lim_{t\rightarrow 0} det \left( A_{\sigma_{1}}(t) A_{\sigma_{2}}(t)^{-1} \right) \\
& = \lim_{t\rightarrow 0} \frac{det\left( \frac{1}{2} \left( I - C(0) \right) + \frac{t}{2} \left( I - \sigma_{1} \right) \right)}
{det\left( \frac{1}{2} \left( I - C(0) \right) + \frac{t}{2} \left( I - \sigma_{2} \right) \right)} \\
& = \lim_{t\rightarrow 0} \frac{det \left( (1+t) I - \left( C(0) + t \sigma_{1} \right) \right) }
{det \left( (1+t) I - \left( C(0) + t \sigma_{2} \right) \right)} \\
& = \lim_{t\rightarrow 0} \frac{det \left( I - C(0) + \frac{t}{t+1} \left( C(0) - \sigma_{1} \right) \right) }
{det \left( I - C(0) + \frac{t}{t+1} \left( C(0) - \sigma_{2} \right) \right)} \\
& = \frac{det \left( C(0) - \sigma_{1} \right) }{det \left( C(0) - \sigma_{2} \right)} \cdot
\lim_{t\rightarrow 0} \frac{det \left( tI + \left( C(0) - \sigma_{1} \right)^{-1} \left( I - C(0) \right) \right) }
{det \left( tI + \left( C(0) - \sigma_{2} \right)^{-1} \left( I - C(0) \right) \right)}.
\endsplit  \tag4.17
$$

\noindent
Using Lemma 4.3  one can check that for $0 < t < 1$
$$ \multline
\left( tI + \left( C(0) - \sigma_{2} \right)^{-1} \left( I - C(0) \right) \right)^{-1} = \\
\left( \frac{1}{t-1} I - \frac{1}{t(t-1)} \left( I + C(0) \right) \left( C(0) - \sigma_{2} \right)^{-1} \right).
\endmultline   $$
Using the above equality we have
$$
\lim_{t\rightarrow 0} \frac{det \left( tI + \left( C(0) - \sigma_{1} \right)^{-1} \left( I - C(0) \right) \right) }
{det \left( tI + \left( C(0) - \sigma_{2} \right)^{-1} \left( I - C(0) \right) \right)} $$
$$ = det \left\{ \left( I + C(0) \right) \left( C(0) - \sigma_{2} \right)^{-1} - \left( C(0) - \sigma_{1} \right)^{-1}
\left( I - C(0) \right) \right\} \qquad $$
$$
= det \left\{ I + \left( C(0) - \sigma_{2} \right)^{-1} \left( C(0) - \sigma_{1} \right)^{-1} \left( I - \sigma_{1}\right)
\left( I + \sigma_{2} \right) \right\}. \tag4.18 $$
Putting
$$\beta = \left( C(0) - \sigma_{2} \right)^{-1} \left( C(0) - \sigma_{1} \right)^{-1} \left( I - \sigma_{1}\right)
\left( I + \sigma_{2} \right) ,$$
we have
$$ \beta = \left( C(0) - \sigma_{2} \right)^{-1} \left( C(0) - \sigma_{1} \right)^{-1} \left( - \sigma_{1} \right) \left( I - \sigma_{1}\right)
\left( I + \sigma_{2} \right). $$
Lemma 4.3 implies that
$$ \beta = -\sigma_{2} \beta = \beta \sigma_{2}. \tag4.19 $$
If we take a basis for $ker B$ by $(+1)$-eigensections and $(-1)$-eigensections of $\sigma_{2}$, $\sigma_{2}$ is of the form
$\left( \smallmatrix I & 0 \\ 0 & - I \endsmallmatrix \right) $. The equation (4.19) implies that $\beta$ has the form
$\left( \smallmatrix 0 & 0 \\ \beta_{3} & 0 \endsmallmatrix \right) $, which shows that
$$ det \left( I + \beta \right) = 1. \tag4.20 $$
From (4.17) and (4.20) we obtain
$$\frac{Det D^{2}_{M, P_{\sigma_{1}}}}{Det D^{2}_{M, P_{\sigma_{2}}}} =
\frac{det \left( C(0) - \sigma_{1} \right) }{det \left( C(0) - \sigma_{2} \right)}, \tag4.21 $$
which completes the proof of the main theorem.

\vskip 1 true cm

\S 5 {\bf Example of Theorem 1.3 on a cylinder }

\vskip 0.5 true cm

In this section we give an example of Theorem 1.3 on a cylinder $[0, r] \times Y$.
Suppose that $\{ f_{1}, \cdots, f_{l}, Gf_{1}, \cdots, Gf_{l} \}$ is an orthonormal basis for $ker B$ and define
a linear map
$\tau : ker B \rightarrow ker B $ by
$$ \tau(f_{i}) = f_{i} , \quad \tau(G f_{i}) = - G f_{i} .$$
For $\theta = (\theta_{1}, \cdots, \theta_{l})$ with $0 < \theta_{i} < \frac{\pi}{2}$ define
$\sigma_{\theta} : ker B \rightarrow ker B $ as follows.
For each $i$,
$$
\sigma_{\theta} ( cos \theta_{i} f_{i} + sin \theta_{i} G f_{i} ) = cos \theta_{i} f_{i} + sin \theta_{i} G f_{i}, $$
$$
\sigma_{\theta} ( - sin \theta_{i} f_{i} + cos \theta_{i} G f_{i} ) = sin \theta_{i} f_{i} - cos \theta_{i} G f_{i}.
$$
Then $\tau$ and $\sigma_{\theta}$ are unitary involutions on $ker B$ anticommuting with $G$.
We use these involutions to define the boundary conditions $ P_{\tau} : = P_{>} + \frac{I - \tau}{2}$ on $Y_{0}$ and
$ P_{\theta} : = P_{<} + \frac{I - \sigma_{\theta}}{2}$ on $Y_{r}$.
We now consider the Dirac Laplacian $(-\partial_{u}^{2}+B^{2})_{P_{\tau}, P_{\theta}}$, whose spectrum is given as follows.
\proclaim{Lemma 5.1}
$$
Spec \left((-\partial_{u}^{2}+B^{2})_{P_{\tau}, P_{\theta}}\right)
 = \{ \mu_{\lambda, j} \mid \lambda\in Spec(B), \lambda \neq 0 , \mu_{\lambda, j} > \lambda^{2} \} $$
$$ \cup \left\{ \left(\frac{\theta_{i} + k \pi}{r}\right)^{2} \mid 1 \leq i \leq l, k = 0, 1, 2, \cdots \right\} $$
$$\cup \left\{ \left(\frac{\pi - \theta_{i} + k \pi}{r}\right)^{2} \mid 1 \leq i \leq l, k = 0, 1, 2, \cdots \right\},
$$
where for $\lambda \neq 0$   $\mu_{\lambda,j}$'s are solutions of
$$\sqrt{\mu - \lambda^{2}} cos (\sqrt{\mu - \lambda^{2}}r) + |\lambda| sin (\sqrt{\mu - \lambda^{2}}r) = 0.$$
\endproclaim

Let $\phi = (\phi_{1},\phi_{2}, \cdots, \phi_{l})$ be another set of angles with $0 < \phi_{i} < \frac{\pi}{2}$.
Then
$$\zeta_{(-\partial_{u}^{2}+B^{2})_{P_{\tau}, P_{\theta}}} (s) - \zeta_{(-\partial_{u}^{2}+B^{2})_{P_{\tau}, P_{\phi}}} (s) $$
$$ = \left( \frac{\pi}{r}\right)^{-2s} \sum_{i=1}^{l}\left\{ \sum_{k=0}^{\infty} \left(k + \frac{\theta_{i}}{\pi}\right)^{-2s}
+ \sum_{k=0}^{\infty} \left(k + \frac{\pi - \theta_{i}}{\pi}\right)^{-2s} \right\} $$
$$ - \left( \frac{\pi}{r}\right)^{-2s} \sum_{i=1}^{l}\left\{ \sum_{k=0}^{\infty} \left(k + \frac{\phi_{i}}{\pi}\right)^{-2s}
+ \sum_{k=0}^{\infty} \left(k + \frac{\pi - \phi_{i}}{\pi}\right)^{-2s} \right\}. $$
The following facts about the Hurwitz zeta-function are well-known ({\it cf.} [AAR] or [MOS]).
\proclaim{Theorem 5.2}
We consider the Hurwitz zeta-function
$$\zeta(s, \alpha) = \sum_{k=0}^{\infty} (k + \alpha )^{-s} \quad \text{ for } \alpha > 0.$$
Then
$$ \zeta(0, \alpha) = \frac{1}{2} - \alpha  \quad \text{ and }
\frac{d}{ds}\zeta(s, \alpha)|_{s=0} = \log \Gamma(\alpha) - \frac{1}{2}\log 2\pi . $$
\endproclaim

\vskip 0.3 true cm

\noindent
Combining Theorem 5.2 with the well-known functional equation
$$\Gamma(z) \Gamma(1-z) = \frac{\pi}{sin \pi z}$$
gives the following result.
\proclaim{Proposition 5.3}
$$ \frac{Det (-\partial_{u}^{2}+B^{2})_{P_{\tau}, P_{\theta}}}{Det(-\partial_{u}^{2}+B^{2})_{P_{\tau}, P_{\phi}}}
= \prod_{i=1}^{l} \frac{sin^{2} \theta_{i}}{sin^{2} \phi_{i}}.$$
\endproclaim

\vskip 0.3 true cm

On the other hand, the scattering matrix $C(\lambda)$ of $G(\partial_{u} + B)_{P_{\tau}}$ on $[0, \infty) \times Y$ is given by
$C(\lambda) = \tau$ and the generalized eigensection $E(f, \lambda)$ is given by
$$
E(f, \lambda) = e^{-i \lambda u} ( f - i G f ) + e^{i \lambda u} \tau ( f - i G f ). $$
Setting
$$F_{i} = cos\theta_{i} f_{i} + sin\theta_{i} G f_{i} \quad  \text{ and } \quad
GF_{i} = - sin\theta_{i} f_{i} + cos\theta_{i} G f_{i}, $$
we have
$$ f_{i} = cos\theta_{i} F_{i} - sin\theta_{i} G F_{i} \quad  \text{ and }  \quad
Gf_{i} = sin\theta_{i} F_{i} + cos\theta_{i} GF_{i} .$$
The ordered basis $\{ f_{1}, Gf_{1}, f_{2}, G f_{2}, \cdots, f_{l}, Gf_{l} \}$ gives the following matrices.
$$
\tau = \left( \matrix
E_{1} & 0 & 0 \\
0 & \ddots & 0 \\
0 & 0 & E_{l} \endmatrix \right)
\text{ with }
E_{i} = \left( \matrix
1 & 0 \\
0 & - 1 \endmatrix \right).  \qquad \qquad
$$

$$
\sigma_{\theta} = \left( \matrix
A_{1} & 0 & 0 \\
0 & \ddots & 0 \\
0 & 0 & A_{l} \endmatrix \right)
\text{ with }
A_{i} = \left( \matrix
cos 2 \theta_{i} & sin 2 \theta_{i} \\
sin 2 \theta_{i} & - cos 2 \theta_{i} \endmatrix \right).
$$
Similarly,
$$
\sigma_{\phi} = \left( \matrix
B_{1} & 0 & 0 \\
0 & \ddots & 0 \\
0 & 0 & B_{l} \endmatrix \right)
\text{ with }
B_{i} = \left( \matrix
cos 2 \phi_{i} & sin 2 \phi_{i} \\
sin 2 \phi_{i} & - cos 2 \phi_{i} \endmatrix \right).
$$

\noindent
Elementary computation shows the following fact.
\proclaim{Proposition 5.4}
$$
\frac{det \left( \tau - \sigma_{\theta} \right)}{det \left( \tau - \sigma_{\phi} \right)}
= \prod_{i=1}^{l} \frac{sin^{2}\theta_{i}}{sin^{2}\phi_{i}}.$$
\endproclaim
\noindent
The equality of Proposition 5.4 and 5.5 gives an example of Theorem 1.3 on a cylinder.

\vskip 1 true cm

\vskip 0.3 true cm
E-mail address : ywlee\@math.inha.ac.kr


\Refs
\widestnumber
\key{APSS}

\ref
\key{AAR}
\by G. E. Andrews, R. Askey and R. Roy
\book Special functions
\publ Cambridge Univ. Press     \yr 1999
\endref



\ref
\key{BFK}
\by D. Burghelea, L. Friedlander and T. Kappeler
\paper Mayer-Vietoris type formula for determinants of elliptic differential operators
\jour J. of Funct. Anal.  \vol 107     \yr 1992    \page 34-66
\endref


\ref
\key{C}
\by J. B. Conway
\book A Course in Funcrional Analysis
\publ Springer-Verlag New York Inc.     \yr 1985
\endref


\ref
\key{KV}
\by M. Kontsevich and S. Vishik
\paper  Geometry of determinants of elliptic operators
\jour Functional analysis on the eve of the 21st century, Vol 1,  Progr. Math.  \vol 131   \yr 1993  \page 173-197
\endref



\ref
\key{L1}
\by Y. Lee
\paper  Mayer-Vietoris formula for the determinants of elliptic operators of Laplace-Beltrami type
(after Burghelea, Friedlander and Kappeler)
\jour Diff. Geom. and Its Appl. \vol 7 \yr 1997 \page 325-340
\endref


\ref
\key{L2}
\by Y. Lee
\paper Burghelea-Friedlander-Kappeler's gluing formula for the zeta determinant and its applications
to the adiabatic decompositions of the zeta-determinant and the analytic torsion
\jour Trans. Amer. Math. Soc. \vol 355-10    \yr 2003  \page 4093-4110
\endref


\ref
\key{L3}
\by Y. Lee
\paper Burghelea-Friedlander-Kappeler's gluing formula and
the adiabatic decomposition of the zeta-determinant of a Dirac Laplacian
\jour Manuscripta Math.   \vol 111    \yr 2003  \page 241-259
\endref



\ref
\key{LW}
\by M. Lesch and K. Wojciechowski
\paper  On the $\eta$-invariant of generalized Atiyah-Patodi-Singer problems
\jour Illinoise J. Math.  \vol 40   \yr 1996  \page 30 - 46
\endref



\ref
\key{M}
\by W. M\"uller
\paper  Eta invariant and manifolds with boundary
\jour J. of Diff. Geom.  \vol 40  \yr 1994 \page 311-377
\endref


\ref
\key{MOS}
\by W. Magnus, F. Oberhettinger and R. P. Soni
\book Formulas and Theorems for the Special functions of Mathematical Physics
\publ Springer-Verlag Berlin Heidelberg New York      \yr 1966
\endref


\ref
\key{PW1}
\by P. Park and K. Wojciechowski
\paper  Analytic surgery of the $\zeta$-determinant of the Dirac operator
\jour Nuclear Phys. B Proc. Suppl.  \vol 104
\yr2002  \page 89-115
\endref



\ref
\key{PW2}
\by P. Park and K. Wojciechowski
\paper  Adiabatic decomposition of the $\zeta$-determinant and scattering theory
\yr2002 \jour preprint
\endref



\endRefs
\enddocument